\documentstyle[11pt]{article}
\begin{document}
\raggedbottom
\setlength{\textwidth}{6.0in}
\setlength{\textheight}{8.5in}
\setlength{\oddsidemargin}{0.2in}
\setlength{\topmargin}{-.05in}

\title{NONSTANDARD TRANSFINITE DIGRAPHS}
\author{A.H. Zemanian}
\date{}

\newcommand{\be}{\begin{equation}}
\newcommand{\ee}{\end{equation}}
\newcommand{\la}{\leftarrow}
\newcommand{\ra}{\rightarrow}
\newcommand{\hla}{hookleftarrow}
\newcommand{\hra}{\hookrightarrow}
\newcommand{\dv}{\dashv}
\newcommand{\vd}{\vdash}
\newcommand{\as}{\asymp}
\newcommand{\lla}{\langle}
\newcommand{\rra}{\rangle}
\newcommand{\N}{I \kern -4.5pt N}
\newcommand{\sss}{^{*}\!}

\maketitle \baselineskip21pt

{\ Abstract --- Nonstandard digraphs and transfinite digraphs have
been defined and examined in two prior technical reports.  The
present work examines digraphs that are both nonstandard and
transfinite.  This requires a combination in certain ways of the
techniques used in the prior two works.  We first construct herein
nonstandard transfinite digraphs for all the natural-number
ranks.  Then, a special kind of nonstandard transfinite digraph
having the "arrow rank" $\vec{\omega}$ needs to be constructed.
Once this is done, the first limit-ordinal rank $\omega$ can be
attained.  This procedure can be continued on to still higher ranks.

}

\section{Introduction}  So far, we have discussed nonstandard generalizations of ordinary digraphs in \cite{nd}
and also transfinite generalizations of ordinary digraphs in \cite{td}.  Our purpose now is to combine these
ideas to generate nonstandard generalizations of transfinite digraphs.   As in \cite{td} ordinary digraphs
will be called 0-digraphs in order to distinguish them from the transfinite digraphs of ranks $\nu$ $(\nu>0$).
We have used an unconventional definition of 0-digraphs \cite[section 2]{nd}, which nonetheless is virtually
equivalent to the usual definition, the only difference being that parallel arcs and arcs which form loops
become allowed.  Our present definition of 0-digraphs conforms with how the $\nu$-digraphs $(\nu>0)$ are set up.

Our notations and symbols herein conform with those used in
\cite{nd} and \cite{td}.  As always, we choose and fix upon a
nonprincipal ultrafilter $\cal F$ on $\N$.

\section{Nonstandard 1-digraphs}

For an ultrapower construction of a nonstandard 1-digraph, we start
with a sequence $\lla D_{n}^{1}: n\in\N\rra$ of (standard)
1-digraphs $D_{n}^{1}$ as defined and analyzed in \cite[Section
3]{td}.  Specifically, for each $n\in\N$, \be
D_{n}^{1}\;=\;\{A_{n},V_{n}^{0},V_{n}^{1}\}  \label{2.1} \ee where
$A_{n}$ is the set of arcs, $V_{n}^{0}$ is the set of 0-vertices,
and $V_{n}^{1}$ is the set of 1-vertices.  Correspondingly, we have
a sequence $\lla D_{n}^{0}\rra$ of 0-digraphs
$D_{n}^{0}=\{A_{n},V_{n}^{0}\}$.  Now, let $T_{n}^{0}$ be the set of
all ditips of $D_{n}^{0}$.  Each ditip of $D_{n}^{0}$ is either an
intip or an outtip.  We refer to these now as 0-{\em ditips}, 0-{\em
intips}, and 0-{\em outtips}, respectively.  It is from these
0-ditips that the 1-vertices in $V_{n}^{1}$ were constructed.  From
the sequence $\lla D_{n}^{0}\rra$ we construct the nonstandard
0-digraph $\sss D^{0}=\{\sss A,\sss V^{0}\}$ exactly as is done in
\cite[Section 3]{nd}.  Here, $\sss A$ is the set of nonstandard arcs
and $\sss V^{0}$ is the set of nonstandard 0-vertices.

Our next objective is to obtain the set $\sss V^{1}$ of "nonstandard
1-vertices."  In accordance with the 1-vertices in $V_{n}^{1}$, we
have $T_{n}^{0}=\cup_{k\in K} T_{n,k}^{0}$, where $K$ is the index
set of the partition of $T_{n}^{0}$ that yielded the 1-vertices
$v_{n,k}^{1}= T_{n,k}^{0}$.  If two 0-ditips $e_{n}$ and $f_{n}$ of
$D_{n}^{0}$ are in the same 1-vertex $v_{n,k}^{1}$ of $D_{n}^{1}$,
we say that $e_{n}$ and $f_{n}$ are {\em shorted together} by
$v_{n,k}^{1}$, and we write $e_{n}\as f_{n}$;  otherwise, if there
is no such 1-vertex, we say that $e_{n}$ and $f_{n}$ are {\em not
shorted together}, and we write $e_{n}\not\as f_{n}$.

Our next step is to make an ultrapower construction to get the
nonstandard 1-vertices.  Consider sequences of 0-ditips such as
$\lla e_{n}\rra$, where $e_{n}$ is a 0-ditip of $D_{n}^{0}$.  Two
such sequences  $\lla e_{n}\rra$ and $\lla f_{n}\rra$ are taken to
be equivalent if $e_{n}=f_{n}$ for almost all $n$.  This partitions
the set of all such sequences into equivalence classes, reflexivity
and symmetry being obvious and transitivity easily shown (i.e., if
$e_{n}=f_{n}$ a.e. and $f_{n}=g_{n}$ a.e., then $e_{n}=g_{n}$ a.e.).
Each such equivalence class is a {\em nonstandard 0-ditip}, and we
denote it by $[e_{n}]$ where $\lla e_{n}\rra$ is a representative
(i.e., is any sequence) in the equivalence class.

Now each (standard) 0-ditip of $D_{n}^{0}$ is either a 0-intip or a
0-outtip.  Given any sequence $\lla e_{n}\rra$ of 0-ditips of
$D_{n}^{0}$, let $N_{i}=\{n: e_{n}$ is a 0-intip$\}$ and $N_{o}=\{n:
e_{n}$ is a 0-outtip$\}$.  Thus, $N_{i}\cap N_{o} =\emptyset$ and
$N_{i}\cup N_{o}=\N$.  So, exactly one of $N_{i}$ and $N_{o}$ is a
member of $\cal F$.  If it is $N_{i}$ (resp. $N_{o}$), $\lla
e_{n}\rra$ is a representative of a {\em nonstandard 0-intip} (resp.
a {\em nonstandard 0-outtip}). As stated in the preceding paragraph,
we refer to both of these as a {\em nonstandard 0-ditip}; now, we
use boldface notation to denote it:  ${\bf e}=[e_{n}]$.

Next, let ${\bf e}=[e_{n}]$ and ${\bf f}=[f_{n}]$ be two nonstandard
0-ditips.  Let $N_{ef}=\{n:e_{n}\as f_{n}\}$ and
$N_{ef}^{c}=\{n:e_{n}\not\as f_{n}\}$.  Exactly one of $N_{ef}$ and
$N_{ef}^{c}$ is a member of $\cal F$.  If it is $N_{ef}$ (resp.
$N_{ef}^{c}$), we say that $\bf e$ and $\bf f$ are {\em shorted
together} (resp. are {\em not shorted together}), and we write ${\bf
e}\as{\bf f}$ (resp. ${\bf e}\not\as {\bf f}$).  Furthermore, we
take it that $\bf e$ is shorted to itself: ${\bf e}\as{\bf e}$.  This
shorting is an equivalence relation for the set of all nonstandard
0-ditips, reflexivity and symmetry being obvious and transitivity
being a result of $\{n: e_{n}\as f_{n}\}\cap\{n:f_{n}\as
g_{n}\}\subseteq \{n: e_{n}\as g_{n}\}$.  The resulting equivalence
classes are the {\em nonstandard 1-vertices}.

This definition can be shown to be independent of the representative
sequences chosen for the nonstandard 0-ditips.  Indeed, let ${\bf
e}=[e_{n}]=[\tilde{e}_{n}]$ and let ${\bf
f}=[f_{n}]=[\tilde{f}_{n}]$.  Set
$N_{e}=\{n:e_{n}=\tilde{e}_{n}\}\in{\cal F}$ and $N_{f}=\{n: f_{n}=\tilde{f}\}\in {\cal F}$.
Assume ${\bf e}\as{\bf f}$; that is $[e_{n}]\as[f_{n}]$.  Thus, $N_{ef}=\{n:
e_{n}\as f_{n}\}\in{\cal F}$.  We want to show that
$N_{\tilde{e}\tilde{f}}=\{n:\tilde{e}_{n}\as \tilde{f}_{n}\}$ is a
member of $\cal F$, which will mean that $[\tilde{e}_{n}]\as
[\tilde{f}_{n}]$, that is, $\tilde{\bf e}\as\tilde{\bf f}$.  We have
$(N_{e}\cap N_{f}\cap N_{ef})\subseteq N_{\tilde{e}\tilde{f}}$, which
implies our desired conclusion.

Altogether, we have defined a nonstandard 1-vertex ${\bf v}^{1}$ to
be any set in the partition of the set of nonstandard 0-ditips
induced by the shorting $\as$.  $\sss V^{1}$ will denote the set of
nonstandard 1-vertices.  Moreover, we have the set $\sss A$ of
nonstandard arcs and the set $\sss V^{0}$ of nonstandard 0-vertices
obtained from the sequence $\lla D_{n}^{0}\rra$, where $D_{n}^{0}=\{ A_{n},V_{n}^{0}\}$,
as above.\footnote{See \cite[Section 3]{td} for a detailed
exposition of this.}  With all these results in hand, we now have
the nonstandard 1-digraph $\sss D^{1}$ as the triplet: \be \sss
D^{1}\;=\;\{\sss A,\sss V^{0},\sss V^{1}\}   \label{2.2} \ee

Let us mention another fact.  We constructed the nonstandard
1-digraph $\sss D^{1}$ by starting with a particular sequence $\lla
D_{n}^{1}\rra$, where each $D_{n}^{1}$ is given by (\ref{2.1}).  But, any
other sequence $\lla \tilde{D}_{n}^{1}\rra$, where
$\tilde{D}_{n}^{1}=\{\tilde{A}_{n},\tilde{V}_{n}^{0},\tilde{V}_{n}^{1}\}$,
would yield the same $\sss D^{1}$ so long as $\{n:
A_{n}=\tilde{A}_{n}\}$, $\{n: V_{n}^{0}=\tilde{V}_{n}^{0}\}$, and
$\{n: V_{n}^{1}=\tilde{V}_{n}^{1}\}$ are all members of $\cal F$.
In this regard, see \cite[Theorem 12.1.1]{go}.

The "underlying nonstandard 1-graph" $\sss G^{1}=\{\sss B, \sss
X^{0},\sss X^{1}\}$ is obtained as follows.  We again start with the
sequence $\lla D_{n}^{1}\rra$, where
$D_{n}^{1}=\{A_{n},V_{n}^{0},V_{n}^{1}\}$.  For each $n$, we remove
the directions of the arcs in $A_{n}$  to get branches.  Two
oppositely directed arcs incident to the same two 0-vertices will
become a pair of parallel branches;  parallel branches are allowed.
In this way, the set $B_{n}$ is obtained from the set $A_{n}$ of
arcs.  Each branch is a pair of $(-1)$-tips, each $(-1)$-tip having
no direction, and we can proceed to get the set $X_{n}^{0}$ of
0-nodes and the set $X_{n}^{1}$ of 1-nodes by using the partitioning
and shorting corresponding to those used to get the set $V_{n}^{0}$
of 0-vertices and the set $V_{n}^{1}$ of 1-vertices.  Then, we can
mimic the construction of the nonstandard arcs and the nonstandard
0-vertices (see \cite[Section 3]{nd}) to get the nonstandard
branches and the nonstandard 0-nodes using the same ultrapower
construction  (i.e., the same partitioning and shortings) as that
used for $\sss D^{0}=\{\sss A,\sss V^{0}\}$.  This yields the set $\sss
B$ of nonstandard branches and the set $\sss X^{0}$ of nonstandard
0-nodes.  Next, by mimicking the construction of nonstandard
1-vertices, we get the nonstandard 1-nodes; these comprise the set
$\sss X^{1}$.  Altogether, we have the {\em underlying nonstandard
1-graph }
\[ \sss G^{1}\;=\;\{\sss B,\sss X^{0},\sss X^{1}\} \]
for the nonstandard 1-digraph $\sss D^{1}$.

\section{Nonstandard $\mu$-digraphs}

We have already constructed a nonstandard 0-digraph $\sss D^{0}$ from
a sequence $\lla D_{n}^{0}:n\in\N\}$ of standard 0-digraphs
$D_{n}^{0}$;  see \cite[Section 3]{nd}.  We have also constructed
in the preceding Section a nonstandard 1-digraph $\sss D^{1}$ from a
sequence $\lla D_{n}^{1}:n\in\N\rra$ of standard 1-digraphs.
Moreover, we have noted that $\sss D^{0}$ does not depend upon the
choice of the sequence $\lla D_{n}^{0}\rra=\lla \{
A_{n},V_{n}^{0}\}\rra$ so long as $A_{n}$ and $V_{n}^{0}$are
specified for almost all $n$, and similarly for $\sss D^{1}$ so long
as $A_{n}$, $V_{n}^{0}$, and $V_{n}^{1}$ are specified for almost
all $n$.

So, now, we shall start with a sequence $\lla
D_{n}^{\mu};n\in\N\rra$ of $\mu$-digraphs, where \be
D_{n}^{\mu}\;=\;\{A_{n},V_{n}^{0},\ldots,V_{n}^{\mu}\}. \label{3.1}
\ee We assume that the
\[ D_{n}^{\mu-1}\;=\;\{ A_{n},V_{n}^{0},\ldots,V_{n}^{\mu-1}\} \]
have already been used to construct the nonstandard
$(\mu-1)$-digraph \be \sss D^{\mu-1}\;=\;\{\sss A,\sss V^{0},\ldots,\sss
V^{\mu-1}\}  \label {3.2} \ee This has already been done for
$\mu-1=0$.  We will now use the sequence $\lla
D_{n}^{\mu}:n\in\N\rra$ to construct the nonstandard $\mu$-digraph
\[ \sss D^{\mu}\;=\;\{\sss A,\sss V^{0},\ldots, \sss V^{\mu-1}, \sss V^{\mu}\} \]
Again, this result will not depend upon which sequence $\lla
D_{n}^{\mu}:n\in\N\rra$ is used so long as the $A_{n}$, $V_{n}^{0}$,
$\dots$ , $V_{n}^{\mu}$ are specified for almost all n.

Let $T_{n}^{\mu-1}$ be the set of all $(\mu-1)$-ditips (i.e., all
$(\mu-1)$-intips and/or all $(\mu-1)$-outtips) in $D_{n}^{\mu-1}$.
These $T_{n}^{\mu-1}$ are not empty because the $D_{n}^{\mu}$ are assumed
to exist.  $T_{n}^{\mu-1}$ has already been partitioned to get the
set of all $\mu$-vertices in $D_{n}^{\mu}$.  Specifically, we have
the partition $T_{n}^{\mu-1}=\cup_{k\in K}T_{n,k}^{\mu-1}$, where
$K$ is the index set for the partition.  Each set of the partition
is a $\mu$-vertex $v_{k,n}^{\mu}=T_{n,k}^{\mu-1}$ in $D_{n}^{\mu}$.
If two $(\mu-1)$-ditips $e_{n}$ and $f_{n}$ of $D_{n}^{\mu-1}$ are
in the same $\mu$-vertex $v_{k,n}^{\mu}$ of $D_{n}^{\mu}$, we will say
that $e_{n}$ and $f_{n}$ are {\em shorted together} by
$v_{k,n}^{\mu}$, and we will write $e_{n}\as f_{n}$;  otherwise, if
there is no such $v_{k,n}^{\mu}$,we will say that $e_{n}$ and
$f_{n}$ are {\em not shorted together}, and we will write
$e_{n}\not\as f_{n}$.

Our eventual aim is to build the nonstandard $\mu$-vertices by means
of an ultrapower construction.  To this end, first consider all
sequences of $(\mu-1)$-ditips of $D_{n}^{\mu-1}$.  Let $\lla
e_{n}\rra$ and $\lla f_{n}\rra$ be two such sequences.  They are
taken to be equivalent if $e_{n}=f_{n}$ for almost all $n$.  This
partitions that set of all such sequences into equivalence
classes---reflexivity, symmetry, and transitivity being clear.  Each
such equivalence class is a {\em nonstandard} $(\mu-1)$-ditip;  we
denote it by $[e_{n}]$, where $\lla e_{n}\rra$ is any sequence in
the equivalence class.  We call $\lla e_{n}\rra$ a representative
sequence for that nonstandard $(\mu-1)$-ditip.

Now, each $(\mu-1)$-ditip of $D_{n}^{\mu-1}$ is either a
$(\mu-1)$-intip or a $(\mu-1)$-outtip.  Given any sequence $\lla
e_{n}\rra$ of $(\mu-1)$-ditips of $D_{n}^{\mu-1}$, let
$N_{i}^{\mu-1}=\{n: e_{n}$ is a $(\mu-1)$-intip$\}$ and
$N_{o}^{\mu-1}=\{n:e_{n}$ is a $(\mu-1)$-outtip$\}$.  Thus,
$N_{i}^{\mu-1}\cap N_{o}^{\mu-1}=\emptyset$  and $N_{i}^{\mu-1}\cup
N_{o}^{\mu-1}=\N$.  So, exactly one of $N_{i}^{\mu-1}$ and
$N_{o}^{\mu-1}$ is a member of $\cal F$.  If it is $N_{i}^{\mu-1}$
(resp. $N_{o}^{\mu-1}$), $\lla e_{n}\rra$ is a representative
sequence of a {\em nonstandard $(\mu-1)$-intip} (resp. a {\em
nonstandard $(\mu-1)$-outtip}).  We call both of these a {\em
nonstandard $(\mu-1)$-ditip} and denote it by ${\bf e}=[e_{n}]$,
where $\lla e_{n}\rra$ is any representative sequence, as above.

We can now obtain the nonstandard $\mu$-vertices as follows.  Let
${\bf e}=[e_{n}]$ and ${\bf f}=[f_{n}]$ be two nonstandard
$(\mu-1)$-ditips.  Let $N_{ef}=\{n: e_{n}\as f_{n}\}$ and
$N_{ef}^{c}=\{n: e_{n}\not\as f_{n}\}$.  Exactly one of $N_{ef}$ and
$N_{ef}^{c}$ is a member of $\cal F$.  If it is $N_{ef}$ (resp.
$N_{ef}^{c}$), we say that $\bf e$ and $\bf f$ are {\em shorted
together} (resp. are {\em not shorted together}, and we write ${\bf
e}\as{\bf f}$ (resp. ${\bf e}\not\as{\bf f}$).  We also take it that
$\bf e$ is {\em shorted to itself}.  It readily follows that $\as$
is an equivalence relation for the set of all nonstandard
$(\mu-1)$-ditips.  The resulting equivalence classes are the {\em
nonstandard $\mu$-vertices}.

Let us check that this definition does not depend upon the chosen
representative sequences for the $(\mu-1)$-ditips:  Assume that the
shorting of $\bf e$ and $\bf f$ was specified with the
representative sequences $\lla e_{n}\rra$ and $\lla f_{n}\rra$,
respectively.  That is, $N_{ef}=\{n: e_{n}\as f_{n}\}\in{\cal F}$.
Let ${\bf e}=[e_{n}]=[\tilde{e}_{n}]$ and ${\bf
f}=[f_{n}]=[\tilde{f}_{n}]$, where $\lla \tilde{e}_{n}\rra$ and
$\lla \tilde{f}_{n}\rra$ are other representative sequences for $\bf e$ and $\bf
f$.  Thus, $N_{e}=\{n: e_{n}=\tilde{e}_{n}\}\in {\cal F}$ and
$N_{f}=\{n: f_{n}=\tilde{f}_{n}\}\in{\cal F}$.  Now, $(N_{e}\cap
N_{f}\cap N_{ef})$ is contained in the set
$N_{\tilde{e}\tilde{f}}=\{n: \tilde{e}_{n}\as\tilde{f}_{n}\}$.
Consequently, $N_{\tilde{e}\tilde{f}}\in{\cal F}$.  This shows that
the shorting ${\bf e}\as{\bf f}$ does not depend upon which
representative sequences are chosen for $\bf e$ and $\bf f$.  A
similar argument shows that if $[e_{n}]\not\as [f_{n}]$, then
$[\tilde{e}_{n}]\not\as[\tilde{f}_{n}]$.

Thus, we now have a valid definition of the nonstandard
$\mu$-vertices as the sets in the chosen partition of the set of all
nonstandard $(\mu-1)$-ditips induced by the shorting $\as$.  $\sss
V^{\mu}$ will denote the set of all $\mu$-vertices.  Appending $\sss
V^{\mu}$ to the sets within the definition (\ref{3.2}) for $\sss
D^{\mu-1}$, we obtain the {\em nonstandard $\mu$-digraph}: \be \sss
D^{\mu}\;=\;\{\sss A,\sss V^{0},\ldots,\sss V^{\mu-1},\sss V^{\mu}\}
\label{3.3} \ee

We obtain the {\em nonstandard underlying $\mu$-graph} $\sss G^{\mu}$
for $\sss D^{\mu}$ as follows:  Again, we start with the sequence
$\lla D_{n}^{\mu}\rra$ (see (\ref{3.1}) above) used for the
construction of $\sss D_{n}^{\mu}$.  For each $n$, we now remove the
directions of all the arcs in $A_{n}$.  This may yield some branches
in parallel, but we allow parallel branches.   $B_{n}$ denotes the
set of all so-obtained branches.  Next, the same shortings that
produced the 0-vertices can again be used to create the 0-nodes, the
set of which we denote by $X_{n}^{0}$.  Furthermore, the same
partitioning and shorting that produced the $k$-vertices
$(k=1,\ldots,n)$ now produces the $k$-nodes, the sets of which are
now denoted by $X_{n}^{k}$.  We can then employ the same
partitionings and shortings used in the ultrapower construction of
the nonstandard $\mu$-digraph $\sss D^{\mu}$ to get the nonstandard
$\mu$-graph $\sss G^{\mu}$.  We say that $\sss G^{\mu}$ {\em
underlies} $\sss D^{\mu}$.

\section{Nonstandard $\vec{\omega}$-digraphs}

As should by now be expected, we start this time with a sequence
$\lla D_{n}^{\vec{\omega}}:n\in\N\rra$, where each
$D_{n}^{\vec{\omega}}$ is an $\vec{\omega}$-digraph having the arrow
rank $\vec{\omega}$.  $\vec{\omega}$-digraphs are discussed in
\cite[Section 6]{td}.  Thus, \be
D_{n}^{\vec{\omega}}\;=\;\{A_{n},V_{n}^{0},V_{n}^{1},\ldots\}
\label{4.1} \ee contains vertex sets $V_{n}^{\mu}$ for all the
natural-number ranks $\mu$.

Corresponding to $D_{n}^{\vec{\omega}}$, we have the $\mu$-digraph
\[ D_{n}^{\mu}\;=\;\{A_{n},V_{n}^{0},V_{n}^{1},\ldots,V_{n}^{\mu}\} \]
for each natural number $\mu$.  Using these, we can construct the
nonstandard $\mu$-digraph $\sss D^{\mu}$ as in the preceding section.
Having done so for every natural number $\mu$, the {\em nonstandard
$\vec{\omega}$-digraph} is simply the infinite set \be \sss
D^{\vec{\omega}}\;=\;\{\sss A,\sss V^{0},\sss V^{1},\ldots\}
\label{4.2} \ee containing $\sss V^{\mu}$ for all the natural numbers
$\mu$.  Here, too, this nonstandard $\vec{\omega}$-digraph will not
depend upon which sequence $\lla D_{n}^{\vec{\omega}}\rra$ is used so
long as the sets $A_{n},V_{n}^{0},V_{n}^{1},\ldots$ remain the same for
almost all $n$.

\section{Nonstandard $\omega$-digraphs}

This time we start with a sequence $\lla D_{n}^{\omega}:n\in\N\rra$
of $\omega$-digraphs \be
D_{n}^{\omega}\;=\;\{A_{n},V_{n}^{0},V_{n}^{1},\dots,V_{n}^{\omega}\}.
\label{5.1} \ee These $\omega$-digraphs are discussed in
\cite[Section 7]{td}.

This implies that we also have the corresponding
$\vec{\omega}$-digraphs
\[ D_{n}^{\vec{\omega}}\;=\;\{A_{n},V_{n}^{0},V_{n}^{1},\ldots\} \]
consisting of all the sets in (\ref{5.1}) except $V_{n}^{\omega}$.
Moreover, we have the $\vec{\omega}$-ditips (i.e., the
$\vec{\omega}$-intips and $\vec{\omega}$-outips) of
$D_{n}^{\vec{\omega}}$; see \cite[Section 6]{td}.

We now proceed exactly
as in Section 3 above, replacing $T_{n}^{\mu-1}$ by the set
$T_{n}^{\vec{\omega}}$, of all $\vec{\omega}$-ditips of
$D_{n}^{\vec{\omega}}$, to get the (standard) $\omega$-vertices in
$D_{n}^{\omega}$ as shortings among the $\vec{\omega}$-ditips.  We
then build the nonstandard $\omega$-vertices by following exactly
the ultrapower construction employed in Section 3.  More specifically, using now
sequences of $\vec{\omega}$-ditips from the $D_{n}^{\vec{\omega}}$
in place of the $(\mu-1)$-ditips, we get the nonstandard
$\vec{\omega}$-ditips, each of which can be identified as either a
nonstandard $\vec{\omega}$-intip or a nonstandard
$\vec{\omega}$-outtip.  Finally, shortings among the nonstandard
$\vec{\omega}$-ditips yield the nonstandard $\omega$-vertices, the
set of which is denoted by $\sss V^{\omega}$.  All this yields the
{\em nonstandard $\omega$-digraph}: \be \sss D^{\omega}\;=\;\{\sss
A,\sss V^{0},\sss V^{1},\ldots,\sss V^{\omega}\}  \label{5.2} \ee

\section{Nonstandard transfinite graphs of still higher ranks}

The same constructions that were used to get the nonstandard
$\mu$-digraphs for all the natural-number ranks $\mu$ and then the
nonstandard digraphs of ranks $\vec{\omega}$ and $\omega$ can now be
used to get the nonstandard digraphs of ranks
$\omega+1,\omega+2,\ldots,\omega+\vec{\omega},\omega\cdot 2$.  This
process can be continued still further to get nonstandard digraphs
of even higher countable-ordinal ranks.  As above, before a
countable limit-ordinal rank, a construction for an arrow rank must
be interposed.

\end{document}